\documentclass[12pt,twoside,leqno]{article}

\usepackage{amsmath}
\usepackage{amssymb}
\usepackage{amsthm}
\usepackage{amscd}
\usepackage{color}
\usepackage{accents}
\usepackage{mathrsfs}
\usepackage{url}

\usepackage{enumitem}
\usepackage{euscript}
\usepackage{multirow}
\usepackage{latexsym}
\usepackage{amsbsy}
\usepackage[pdftex]{graphicx}
\usepackage{epsfig}
\usepackage{subdepth}

\usepackage{float}
\usepackage{color}
\usepackage{relsize}
\usepackage[geometry]{ifsym}
\usepackage{pdfpages}
\usepackage[makeroom]{cancel}
\usepackage[normalem]{ulem}

\usepackage[pagebackref=false,colorlinks=true,linkcolor=blue,citecolor=blue,urlcolor=blue]{hyperref}

\numberwithin{equation}{section}

\newcommand{\Proof}{\noindent\textit{Proof.} }
\def\C{{\mathbb C}}

\def\R{{\mathbb R}}

\def\dd{\,\mathrm{d}}
\newcommand{\diag}{\mathop{\mathrm{diag}}}
\newcommand{\etr}{\mathop{\mathrm{etr}}}
\newcommand{\tr}{\mathop{\mathrm{tr}}}

\def\calS{\mathcal{S}}
\newcommand{\calSd}{{\mathcal{S}^{d \times d}}}
\newcommand{\calSdpd}{{\mathcal{S}_+^{d \times d}}}

\def\fhat{\widehat{f}}
\def\fonehat{\widehat{f_1}}
\def\ftwohat{\widehat{f_2}}
\def\fthreehat{\widehat{f_3}}
\def\fjhat{\widehat{f_j}}

\newtheorem{theorem}{Theorem}[section]

\theoremstyle{definition}

\theoremstyle{remark}

\newtheorem{example}[theorem]{\rm\bf Example}

\begin{document}

\title{\Large\textbf{A Tribute to Richard Askey, and an Expansive View of Some of his Favorite Beta Integrals}}

\author{Donald Richards}

\date{ }

\maketitle

\begin{abstract}
This article represents a personal tribute to Richard Askey together with a new look at some of his favorite integrals, including the Cauchy beta integral.  The article also provides some new multidimensional extensions of Cauchy's beta integral in which the domain of integration is the space of real symmetric matrices, and these multidimensional integrals are used to obtain some special cases of the Cauchy--Selberg integrals.  


\end{abstract}

\maketitle

\markboth{\hfill{\rm\sc{D. Richards}}\hfill}{\hfill{\rm\sc{A Tribute to Richard Askey}}\hfill}

\tableofcontents

\normalsize



\section{Reminiscences}
\label{sec_reminiscences}

In late October, 1982 my friend and co-author, Kenneth Gross, had been scheduled to deliver a lecture at the American Mathematical Society's 797th meeting, to be held at College Park, Maryland.  The lecture was to be on work by the two of us and it was to be given in the special session on Lie groups and generalized classical special functions and on the topic of generalized Bessel functions of matrix argument arising in harmonic analysis, physics, and statistics.  

Then at a late stage, too late for me to be listed on the program as a speaker, Ken suggested that I should give the talk.  The suggestion made me nervous, for I guessed that the audience was likely to know far more about Bessel functions than I did.  In preparing this tribute to Dick Askey I consulted the \textit{Notices} \cite[p.~531]{Notices}, and I see even now that I was right to be unnerved, for the list\footnote{The complete list was: L. Biedenharn, Richard Askey, L. Auslander, Kenneth Baclawaki, Johan G. F. Belinfante, Charles F. Dunkl, Loyal Durand, Jacques Faraut, Philip Feinsilver, Daniel Flath, Bruno Gruber, Kenneth I. Gross, R. A. Gustafson, Robert Hermann, Wayne J. Holman III, F. T. Howard, Yehiel Ilamed, Kenneth D. Johnson, Tom H. Koornwinder, W. D. Lichtenstein, James D. Louck, Willard Miller, Jr., S. C. Milne, Walter G. Morris, Donald Richards, Walter Schempp, Michael F. Singer, Audrey A. Terras, and Jacob Towber.} of organizers, speakers, and authors at the special session included numerous leaders in the field of special functions.  

In any event, I gave my talk and it seemed to go well; I even had a few questions from the audience.  After my talk ended, a gently smiling man approached me and introduced himself as Dick Askey.  By the end of the meeting, Dick had invited me to visit him at Madison, which I did in February, 1983.  I always recall that Dick was the first ever to invite me to lecture at another school, and that he eagerly introduced me to I. J. Schoenberg and to R. C. Buck.  The introduction to Schoenberg turned out eventually to be highly beneficial, for it partly led me years later into research on the theory of total positivity and the theory of positive definite functions \cite{pdfunctions,tp1}, both of which are areas in which Schoenberg played pioneering roles.  

Dick also introduced me to Dennis Stanton who invited me to visit him in Minneapolis.  I received from Dennis superb advice about the literature on root systems and finite reflection groups, and that advice spawned much later the first connection \cite{tp2} between the areas of total positivity and finite reflection groups.

During the 1983 visit to Madison, Dick introduced me to Selberg's integral and urged me to think about possible applications to statistics, which led in due course to the articles \cite{AskeyRichards,Richards_IMA}.  On the other hand, when Dick suggested that I study the area of $q$-series, I politely declined to do so because ``$q$-series will forever be too difficult for me.''  Dick then noted that ``some $q$-hypergeometric series are special cases of the hypergeometric functions of Hermitian matrix argument, so the problems that you and Ken Gross are studying are more difficult than some $q$-series problems.''  I explained that it may well be true that the special functions of matrix argument are more difficult, but there is more ``structure to the theory underpinning the matrix-argument functions,'' whereas ``the theory of $q$-series seems to me to be the mathematical equivalent of close-up magic tricks: I see it, but it defies the imagination, and even after seeing it, I can't do it.''  I realized many years later that Dick, by encouraging me to look into the area of $q$-series, was detecting more promise in me than I was seeing in myself.

On the other hand, Dick did get me to look closely at some of Ramanujan's Notebooks \cite{Berndt}.  That motivated Ding, Gross, and me later to derive a generalization \cite{DGR} of Ramanujan's Master Theorem within the setting of symmetric cones, a result that has spawned numerous papers by others.  

Over the years, Dick told me some hilarious stories about his experiences outside of research.  My favorite such story concerns the aftermath of his trip to attend the 1966 International Congress of Mathematicians (ICM) in Moscow.  While attending the ICM, Dick bought a conspicuous {\em ushanka}-style fur hat 
and a Russian winter coat.  A few years later, Waleed Al-Salam invited Dick to give a talk at the University of Alberta, in Edmonton.  So Dick went off to Edmonton and, in preparation for the frigid winter temperatures awaiting him, he wore his Russian hat and coat.  As there were few direct flights to Edmonton, Dick flew to Calgary and then changed planes for his continuing trip to Edmonton.  After Dick gave his lecture he went that evening to Al-Salam's home for dinner, and then he flew the next day back to the U.S.  By a slim coincidence there happened to be in Calgary at that time a Russian athletics team, and one of the Russian athletes disappeared.  The Canadian police, in a frantic search for the missing athlete, thought that Dick was the man they were looking for; after all, Dick fit the profile well given his eye-catching Russian hat, Russian coat, and his athletic build (Dick was then in his early forties).  So the next day, the Edmonton police knocked on Al-Salam's door and insisted that he tell them where they could find the Russian athlete.  Al-Salam was dumbfounded, naturally, and became even more so as the police continued to claim strenuously that Al-Salam surely knew of the athlete's location.  Eventually Al-Salam realized that the police had mistaken Dick for the missing athlete, so he explained to the police that, first, they had mistakenly confused Dick with the missing athlete and, second, they had also missed Dick's passage through two airport security zones and his departure to the U.S.\footnote{During much laughter between Dick and me, I suggested that he ought to have called the Edmonton police and said: ``You may have found me `spooky,' but it was freezing in Moscow when I bought the hat and coat and I was simply trying to be an `intelligent' guy and come in from the cold.''  More laughter ensued, for we both knew that the police would not be amused by a thinly-veiled reference to Le Carr\'e \cite{LeCarre}.}  

On several occasions, Dick served as my guardian angel.  Despite my hesitation about $q$-series, Dick suggested in early 1985 that I attend an upcoming conference on $q$-series \cite{Andrews1985}, to be held in May, 1985 at Arizona State University.  At that conference I met Joe Gillis, whose work on random walks and combinatorial applications of orthogonal polynomials \cite{Zeilberger} I had admired for years.  Gillis was pleased to learn during a lunch-time chat that I had read his work on random walks, and he kindly described for me what he saw as his more important contributions to mathematics and to science in Israel, and that too under difficult working conditions.  That conversation with Gillis was influential, for it helped me to remember that my problems were far easier than those faced by many other mathematicians, and I tried to keep that in mind as I continued to work on the Selberg integrals \cite{Richards_IMA}.  

During a possibly life-saving conversation with Dick at that 1985 conference, I mentioned to him that I was interested in taking a helicopter tour of the Grand Canyon.  Dick suddenly looked at me pointedly and said, ``don't do that,'' and he went on to explain that too often there was national news of a fatal accident involving aircraft flying over the Grand Canyon.  Naturally, I took his advice; and later I was very glad to have taken his advice as there occurred in 1986 another such accident.  

In October, 2009 Dick attended the American Mathematical Society meeting at University Park, Pennsylvania and my wife, Mercedes, and I invited him to dinner at our home.  During dinner, Dick noticed a wall-mounted plaque that Mercedes received in 2008 when the Institute of Jamaica, which is based in Kingston, Jamaica, awarded her the Institute's Musgrave Gold Medal for her work in astronomy and astrophysics.  Naturally, I told Dick that was impressed by the award.  At that point, Dick chided me gently that I ought to be far, far more impressed than I appeared to be.  Dick told us that there had been a linguist at the University of Wisconsin who had also been awarded the Institute's Musgrave Gold Medal, and Dick then explained that the linguist had been in his mid-70's at the time that he received his award, whereas Mercedes then was in her early 50's.  Mercedes and I, both having been born in Jamaica, were stunned that this guy, Askey, knew more about the Institute of Jamaica than we did!

But let there be no mistake: Dick could be candid when he felt that the occasion called for frankness.  On one occasion Dick disapproved in a blunt opinion \cite{Askey_MR} published in \textit{Mathematical Reviews}, no less, of a series of articles, none of which he felt should have been published.  Also, Dick's widely-disseminated, sharply-worded opinions on mathematics pedagogy \cite{Askeyteaching} attracted much attention \cite{AddRoit}.  As for Dick's comments to me about the continual lowering of standards for collegiate mathematics majors, his concerns have proved to be astute.  At the same time, Dick provided extensive and loving reviews of mathematics books for children \cite{Askeyreviews}, and he was from many accounts unstinting of his time with others whom he recognized as being utterly serious about mathematics.  

When I learnt in August, 2019 that Dick had been admitted to hospice, I immediately contacted his family and made plans to fly to Madison to see him.  (As I commented to my family, the principle of ``regret minimization'' required that I visit Dick immediately.)  I met Dick on the morning of September 9, 2019 and, to my pleasant surprise, he had me wheel him around the facility, all the while introducing me to the staff.  I honestly felt like a new grandchild!

In looking back on 1982--2019, the years that I knew Dick, it is impossible for me to thank him enough.  He was a mentor, an adviser, and a friend, and I will forever owe him my deepest gratitude.  When my wife passed away in 2016 after a long illness, Dick called me immediately on hearing the news, and his kind words of condolence were deeply heartfelt by my daughters and me.  His comment in 1982 that my talk at College Park was well-received by the audience constituted enormous encouragement to a neophyte trying to find his niche in research, and in life too.  And even today, more than 42 years later, I continue to work and publish on the Bessel functions of matrix argument.

\section{An expansive view of Cauchy's beta integral}
\label{sec_Cauchy}

We now consider a variety of integrals of a genre that have become attached to Askey's name.  In all of these integrals, it is understood that fractional powers are assigned their principal values and are continued through continuity.  

Cauchy's beta integral \cite{Cauchy1,Cauchy2}, which is stated immediately below, is the first of several examples of beta integrals that we consider in this article.  

\medskip

\begin{example} 
\label{ex_Cauchy_beta}
(Cauchy \cite{Cauchy1,Cauchy2}) 
\textit{Let $\alpha_1,\alpha_2 \in \C$ where $\Re(\alpha_1+\alpha_2) > 1$, $\Re(\sigma_1) > 0$, and $\Re(\sigma_2) > 0$.  Then 
\begin{equation}
\label{eq_Cauchy_beta_integral}
\int_{-\infty}^\infty (1+i\sigma_1 t)^{-\alpha_1} (1-i\sigma_2 t)^{-\alpha_2} \dd t 
= 2\pi \frac{\Gamma(\alpha_1+\alpha_2-1)}{\Gamma(\alpha_1) \Gamma(\alpha_2)} \, \frac{\sigma_1^{\alpha_2-1} \sigma_2^{\alpha_1-1}}{(\sigma_1+\sigma_2)^{\alpha_1+\alpha_2-1}}
\end{equation}
and 
\begin{equation}
\label{eq_Cauchy_beta_integral2}
\int_{-\infty}^\infty (1+i\sigma_1 t)^{-\alpha_1} (1+i\sigma_2 t)^{-\alpha_2} \dd t = 0.
\end{equation}
}\end{example} 

The formula \eqref{eq_Cauchy_beta_integral} was one of Askey's favorite integrals.  As Rahman and Suslov \cite[p.~652]{RahmanSuslov} commented, ``Askey's name seems to be attached to most of the interesting modern extensions of the beta integral, $\ldots$'' so it is not surprising that Askey's name appears often in the literature on Cauchy's beta integral or its $q$-extensions; see \cite{AR,Askey4,Askey5,Askey1,Askey3,AskeyRoy,Habsieger,Wilson}.  The integral \eqref{eq_Cauchy_beta_integral} has also taken on a life of its own outside the domain of the classical special functions, for it has appeared perhaps unexpectedly in distance correlation theory \cite{DER}, diffusion processes \cite{AristaDemni}, and other areas.

As regards proofs of \eqref{eq_Cauchy_beta_integral} and \eqref{eq_Cauchy_beta_integral2} other than those in \cite{Cauchy1,Cauchy2}, we note that Andrews, Askey, and Roy \cite[p.~48]{AAR} provides a proof that is based on recurrence relations as a function of $\alpha_1$ and $\alpha_2$ for the left-hand side of \eqref{eq_Cauchy_beta_integral}.  It appears that the most widely-used approach currently is to transfer the integral to the complex plane and then to apply the calculus of residues; cf.~Rahman and Suslov \cite{RahmanSuslov}.  

What seems not to be well-known, or at least not to have been widely applied, is the result that both \eqref{eq_Cauchy_beta_integral} and \eqref{eq_Cauchy_beta_integral2} are consequences of Parseval's identity for the Fourier transform.  Although we hesitate to claim that this observation is new, we have not found during an extensive search of the extant literature on beta integrals any applications of Parseval's identity for the Fourier transform.  In particular, a well-known biography \cite{Belhoste} of Cauchy does not appear to indicate that Cauchy himself used an approach based on Parseval's identity.

A close reading of Askey and Roy \cite{AskeyRoy} also reveals that they applied the Parseval identity for the Mellin transform to evaluate Barnes' beta integral and several of its $q$-extensions, so we infer that Askey likely was well aware that the Parseval identity for the Fourier transform could be applied to evaluate numerous beta integrals.  We also remark that Paris and Kaminski \cite[p.~83]{ParisKaminski} states that Marichev \cite{Marichev} applied the Parseval identity for the Mellin transform to evaluate numerous integrals; i.e., Marichev concentrated on the Mellin transform rather than on the Fourier transform.  

For $p \ge 1$ let $L^p(\R)$ denote the space of Lebesgue measurable functions $f:\R \to \C$ such that $\int_\R |f(x)|^p \dd x$ is finite.  The space $L^2(\R)$, endowed with the inner product 
$$
\langle f_1,f_2 \rangle = \int_{-\infty}^\infty f_1(x) \overline{f_2(x)} \dd x,
$$
$f_1, f_2 \in L^2(\R)$, is a complex Hilbert space.  For $f \in L^1(\R) \cap L^2(\R)$, define its Fourier transform, 
$$
\fhat(t) = \int_{-\infty}^\infty e^{-itx} f(x) \dd x.
$$
Also, for Fourier transforms $\fonehat,\ftwohat \in L^2(\R)$, one defines the inner product 
$$
\langle\fonehat,\ftwohat\rangle = \int_{-\infty}^\infty \fonehat(t) \overline{\ftwohat(t)} \dd t.
$$
According to the Parseval identity, which establishes that the Fourier transform is a unitary operator on $L^2(\R)$, one has 
\begin{equation}
\label{eq_Parseval}
\frac{1}{2\pi} \langle\fonehat,\ftwohat\rangle = \langle f_1,f_2 \rangle.
\end{equation}

Throughout this article, we usually begin the investigation of each integral by assuming that several parameters are real.  Once an integral is derived for the case of real parameters then we will apply analytic continuation to extend the domain of validity of the integral to regions of analyticity in the complex plane.  For example, once we have established \eqref{eq_Cauchy_beta_integral} for all $\alpha_1,\alpha_2 > 1/2$, and $\sigma_1,\sigma_2 > 0$ then, on observing that the right-hand side of \eqref{eq_Cauchy_beta_integral} is analytic on the larger region $\Re(\alpha_1+\alpha_2) > 1$, $\Re(\sigma_1) > 0$, and $\Re(\sigma_2) > 0$, it follows that the equality holds throughout the larger region.  

\medskip

\noindent\textit{Proof of Example \ref{ex_Cauchy_beta}.}  To establish \eqref{eq_Cauchy_beta_integral} using the Parseval identity, suppose that $\alpha_1,\alpha_2 > 1/2$, and $\sigma_1,\sigma_2 > 0$.  For $j=1,2$, define 
\begin{equation}
\label{eq_f1_function}
f_j(x) = 
\begin{cases}
\sigma_j^{-\alpha_j} x^{\alpha_j-1} \exp(-\sigma_j^{-1}x)/\Gamma(\alpha_j), & x > 0 \\
0, & \hbox{otherwise}
\end{cases}
\end{equation}
so that each $f_j$ is the familiar probability density function of a gamma-distributed random variable.  An elementary calculation yields 
\begin{equation}
\label{eq_f1_function_FT}
\fjhat(t) = (1+i\sigma_j t)^{-\alpha_j},
\end{equation}
$t \in \R$, $j=1,2$.  Since each $\alpha_j > 1/2$ then $f_j \in L^2(\R)$, hence $\fjhat \in L^2(\R)$, so we may apply Parseval's identity, obtaining 
\begin{align*}
\frac{1}{2\pi} \int_{-\infty}^\infty (1+i\sigma_1 t)^{-\alpha_1} (1-i\sigma_2 & t)^{-\alpha_2} \dd t \\
&= \frac{1}{2\pi} \langle\fonehat,\ftwohat\rangle = \langle f_1,f_2 \rangle \\
&= \frac{\sigma_1^{-\alpha_1} \sigma_2^{-\alpha_2}}{\Gamma(\alpha_1)\Gamma(\alpha_2)} \int_0^\infty x^{\alpha_1+\alpha_2-2} \exp(-(\sigma_1^{-1}+\sigma_2^{-1})x) \dd x \\
&= \frac{\sigma_1^{-\alpha_1} \sigma_2^{-\alpha_2}}{\Gamma(\alpha_1)\Gamma(\alpha_2)} (\sigma_1^{-1}+\sigma_2^{-1})^{-(\alpha_1+\alpha_2-1)} \Gamma(\alpha_1+\alpha_2-1),
\end{align*}
which reduces easily to \eqref{eq_Cauchy_beta_integral}

As regards \eqref{eq_Cauchy_beta_integral2}, one defines 
$$
f_3(x) = \begin{cases}
0, & x \ge 0 \\
f_2(-x), & x < 0
\end{cases},
$$
so that 
$$
\fthreehat(t) = \ftwohat(-t) = (1+i\sigma t)^{-\beta},
$$
$t \in \R$.  Noting that $f_1(x) f_3(x) \equiv 0$ then \eqref{eq_Cauchy_beta_integral2} also follows from Parseval's identity, \eqref{eq_Parseval}.  
$\qed$

\medskip

Parseval's identity enables us to adopt an expansive view of Cauchy's beta integral.  Thus several observations can be made immediately.  First, we see that Cauchy's beta integral extends to cases in which the functions $f_1$ and $f_2$ are of any functional form such that the product $f_1 f_2$ also is of a similar form and such that each $\fjhat$ can be calculated explicitly.  We now provide some examples, and in so doing, for $\Re(\alpha), \Re(\beta) > 0$, let 
$$
B(\alpha,\beta) = \frac{\Gamma(\alpha) \Gamma(\beta)}{\Gamma(\alpha+\beta)}
$$
denote the classical beta function.  

\medskip

\begin{example}
\label{example_beta}
\textit{Let $\alpha_1,\alpha_2,\sigma_1,\sigma_2 \in \C$ be such that $\Re(\alpha_1+\alpha_2) > 1$, $\Re(\sigma_1+\sigma_2) > 1$, $\Re(\alpha_1+\sigma_1) > 0$, and $\Re(\alpha_2+\sigma_2) > 0$.  Then 
\begin{multline}
\label{eq_beta_1F1}
\frac{1}{2\pi} \int_{-\infty}^\infty {}_1F_1(\alpha_1;\alpha_1+\sigma_1;-it) \, {}_1F_1(\alpha_2;\alpha_2+\sigma_2;it) \dd t \\
= \frac{\Gamma(\alpha_1+\alpha_2-1)\Gamma(\sigma_1+\sigma_2-1)\Gamma(\alpha_1+\sigma_1)\Gamma(\alpha_2+\sigma_2)}{\Gamma(\alpha_1+\alpha_2+\sigma_1+\sigma_2-2)\Gamma(\alpha_1)\Gamma(\sigma_1)\Gamma(\alpha_2)\Gamma(\sigma_2)}.
\end{multline}
Further 
\begin{equation}
\label{eq_beta_1F12}
\int_{-\infty}^\infty {}_1F_1(\alpha_1;\alpha_1+\sigma_1;it) \, {}_1F_1(\alpha_2;\alpha_2+\sigma_2;it) \dd t = 0.
\end{equation}
}\end{example}

\Proof  For $j=1,2$, suppose initially that $\alpha_j,\sigma_j \in \R$ and that $\alpha_j,\sigma_j > 1/2$.  Consider the probability density function of a beta-distributed random variable, 
\begin{equation}
\label{eq_f_function}
f_j(x) = \begin{cases}
\dfrac{x^{\alpha_j-1} \, (1-x)^{\sigma_j-1}}{B(\alpha_j,\sigma_j)}, & 0 < x < 1 \\
0, & \hbox{otherwise}
\end{cases}.
\end{equation}
It is straightforward to verify that $f_j \in L^1(\R) \cap L^2(\R)$.  Also, it is well-known \cite[p.~1023]{Gradshteyn} that 
$$
\fjhat(t) = {}_1F_1(\alpha_j;\alpha_j+\sigma_j;-it)
$$
for all $t \in \R$.  On applying Parseval's identity, we obtain 
\begin{align*}
\frac{1}{2\pi} \int_{-\infty}^\infty {}_1F_1(\alpha_1;\alpha_1+\sigma_1;-it) & \, {}_1F_1(\alpha_2;\alpha_2+\sigma_2;it) \dd t \\
&= \frac{1}{2\pi} \langle\fonehat,\ftwohat\rangle = \langle f_1,f_2 \rangle \nonumber \\
&= \frac{1}{B(\alpha_1,\sigma_1) B(\alpha_2,\sigma_2)} \int_0^1 x^{\alpha_1+\alpha_2-2} (1-x)^{\sigma_1+\sigma_2-2} \dd x \\
&= \frac{B(\alpha_1+\alpha_2-1,\sigma_1+\sigma_2-1)}{B(\alpha_1,\sigma_1) B(\alpha_2,\sigma_2)}.
\end{align*}
On writing each of these beta functions in terms of gamma functions, we obtain the right-hand side of \eqref{eq_beta_1F1}, and then it follows by analytic continuation that the final result remains valid for the stated range of the parameters $\alpha_j, \beta_j$, $j=1,2$.  

Finally, we obtain \eqref{eq_beta_1F12} by the same argument that led to \eqref{eq_Cauchy_beta_integral2}.  
$\qed$

\medskip

Note that \eqref{eq_Cauchy_beta_integral} is a limiting case of \eqref{eq_beta_1F1}.  This is obtained by setting $\sigma_1 = \sigma_2$ and substituting $t = \sigma_1 u$ in \eqref{eq_beta_1F1} and then letting $\sigma_1 \to \infty$.   

\medskip

\begin{example}
\label{example_1F1}
\textit{Suppose that $\alpha_1,\alpha_2,\gamma_1,\gamma_2,\sigma_1,\sigma_2,\theta \in \C$ such that $\Re(\alpha_1+\alpha_2) > 1$, $\Re(\sigma_1) > 0$, $\Re(\theta) > \Re(\sigma_2) > 0$, and $-\gamma_2$ is not a nonnegative integer.  Then 
\begin{multline}
\label{eq_2f1}
\frac{1}{2\pi} 
\int_{-\infty}^\infty (1+i\sigma_1 t)^{-\alpha_1} (1-i\sigma_2 t)^{-\alpha_2} \, {}_2F_1\big(\gamma_1,\alpha_2;\gamma_2;\theta^{-1}\sigma_2(1-i\sigma_2 t)^{-1}\big) \dd t \\
= \frac{\Gamma(\alpha_1+\alpha_2-1)}{\Gamma(\alpha_1) \Gamma(\alpha_2)} \frac{\sigma_1^{\alpha_2-1} \sigma_2^{\alpha_1-1}}{(\sigma_1+\sigma_2)^{\alpha_1+\alpha_2-1}} \, {}_2F_1\big(\gamma_1,\alpha_1+\alpha_2-1;\gamma_2;\theta^{-1}(\sigma_1^{-1}+\sigma_2^{-1})^{-1}\big).
\end{multline}
}\end{example}

\Proof 
Suppose that $\alpha_1 > 1/2$, $\alpha_2 > 1/2$, $\sigma_1 > 0$, $\theta > \sigma_2 > 0$, and $\gamma_2 \ge \gamma_1 > 0$.  Let $f_1$ be the function in \eqref{eq_f1_function}, with Fourier transform given in \eqref{eq_f1_function_FT}.  Also define 
\begin{equation}
\label{eq_f2_ex_1F1}
f_2(x) = \begin{cases}
c \: x^{\alpha_2-1} \, \exp(-\sigma_2^{-1}x) \, {}_1F_1(\gamma_1;\gamma_2;\theta^{-1}x), & x > 0 \\
0, & \hbox{otherwise}
\end{cases},
\end{equation}
where the normalizing constant $c$ satisfies 
\begin{align*}
c^{-1} &= \int_0^\infty x^{\alpha_2-1} \exp(-\sigma_2^{-1}x) \, {}_1F_1(\gamma_1;\gamma_2;\theta^{-1}x) \dd x \\
&= \Gamma(\alpha_2) \sigma_2^{\alpha_2} \, {}_2F_1(\gamma_1,\alpha_2;\gamma_2;\theta^{-1}\sigma_2).
\end{align*}
This calculation of $c$ also shows that $f_2 \in L^1(\R)$.  

Next, we show that $f_2 \in L^2(\R)$. Indeed, since $\gamma_2 \ge \gamma_1 > 0$ then it follows from a Taylor--Maclaurin expansion that, for $x \ge 0$, 
$$
{}_1F_1(\gamma_1;\gamma_2;x) = \sum_{j=0}^\infty \frac{(\gamma_1)_j}{(\gamma_2)_j} \frac{x^j}{j!} \le \sum_{j=0}^\infty \frac{x^j}{j!} = \exp(x).
$$
Applying this inequality to \eqref{eq_f2_ex_1F1}, we obtain 
\begin{align*}
\int_0^\infty [f_2(x)]^2 \dd x &\le c^2 \int_0^\infty x^{2\alpha_2-2} \exp(-2\sigma_2^{-1}x) \exp(2\theta^{-1}x) \dd x \\
&= c^2 \int_0^\infty x^{2\alpha_2-2} \exp(-2(\sigma_2^{-1}-\theta^{-1})x) \dd x,
\end{align*}
and this integral is finite since, by hypothesis, $\alpha_2 > 1/2$ and $\sigma_2^{-1}-\theta^{-1} > 0$.  

Further, 
\begin{align*}
\ftwohat(t) &= c \, \int_0^\infty x^{\alpha_2-1} \exp(-\sigma_2^{-1}(1+i\sigma_2 t)x) \, {}_1F_1(\gamma_1;\gamma_2;\theta^{-1}x) \dd x \\
&= c \, \Gamma(\alpha_2) \sigma_2^{\alpha_2} (1+i\sigma_2 t)^{-\alpha_2} \, {}_2F_1(\gamma_1,\alpha_2;\gamma_2;\theta^{-1}\sigma_2(1+i\sigma_2 t)^{-1}),
\end{align*}
with convergence for all $t \in \R$ subject to the stated conditions on $\alpha_2$, $\sigma_2$, $\theta$, $\gamma_1$, and $\gamma_2$.  

We have now shown that $f_j \in L^1(\R) \cap L^2(\R)$, hence $\fjhat \in L^2(\R)$, $j=1,2$.  Therefore, by Parseval's identity, 
\begin{align*}
\frac{c \, \Gamma(\alpha_2) \sigma_2^{\alpha_2}}{2\pi} \int_{-\infty}^\infty & 
(1+i\sigma_1 t)^{-\alpha_1} (1-i\sigma_2 t)^{-\alpha_2} \, {}_2F_1\big(\gamma_1,\alpha_2;\gamma_2;\theta^{-1}\sigma_2(1-i\sigma_2 t)^{-1}\big) \dd t \\
&= \frac{1}{2\pi} \langle\fonehat,\ftwohat\rangle = \langle f_1,f_2 \rangle \\
&= \frac{c \, \sigma_1^{-\alpha_1}}{\Gamma(\alpha_1)} \int_0^\infty x^{\alpha_1+\alpha_2-2} \exp\big(-(\sigma_1^{-1}+\sigma_2^{-1})x\big) \, {}_1F_1(\gamma_1;\gamma_2;\theta^{-1}x) \dd x \\
&= \frac{c \, \sigma_1^{-\alpha_1}}{\Gamma(\alpha_1)} \Gamma(\alpha_1+\alpha_2-1) \, (\sigma_1^{-1}+\sigma_2^{-1})^{-(\alpha_1+\alpha_2-1)} \\
& \qquad\qquad\qquad\cdot \, {}_2F_1\big(\gamma_1,\alpha_1+\alpha_2-1;\gamma_2;\theta^{-1}(\sigma_1^{-1}+\sigma_2^{-1})^{-1}\big),
\end{align*}
the latter equality following from a well-known result for the Laplace transform of a confluent hypergeometric function.  
Simplifying this result, we obtain \eqref{eq_2f1} for real values of the parameters, and by applying analytic continuation we obtain the result for the full domain of analyticity.  
$\qed$

\medskip

We remark that the latter result can also be obtained by expanding the ${}_2F_1$ as a hypergeometric series and applying \eqref{eq_Cauchy_beta_integral} to integrate term-by-term.  More generally, the result \eqref{eq_2f1} can be extended to arbitrary ${}_pF_q$ generalized hypergeometric series.  

\medskip

\begin{example}
\label{ex_t_distn}
{\rm 
Consider the probability density function of a renormalized Student's $t$-distribution with $2\nu$ degrees-of-freedom, 
\begin{equation}
\label{eq_t_pdf}
f_\nu(x) = \dfrac{(1+x^2)^{-(\nu+\frac12)}}{B(\nu,\tfrac12)},
\end{equation}
$x \in \R$, $\nu > 0$.  Denoting by $K_\nu(\cdot)$ the modified Bessel function (of the second kind) of order $\nu$, the Fourier transform of $f_\nu$ is 
\begin{equation}
\label{eq_t_cdf}
\fhat_\nu(t) = \begin{cases}
\dfrac{|t|^\nu K_\nu(|t|)}{2^{\nu-1} \Gamma(\nu)}, & t \neq 0 \\
1, & t=0
\end{cases},
\end{equation}
a result that follows from \textit{Basset's integral}: For $t > 0$, 
$$
K_\nu(t) = \pi^{-1/2} \Gamma(\nu+\tfrac12) (2/t)^\nu \int_0^\infty (1+x^2)^{-(\nu+\frac12)} \cos(tx) \dd x.
$$
We refer to Olver and Maximon \cite[eq.~(10.32.11)]{OlverMaximon} and Watson \cite[p.~172]{Watson} for more details about Basset's integral and to Hurst \cite{Hurst} for an alternative derivation of \eqref{eq_t_cdf}.  

Suppose that $\nu_1,\nu_2 > 0$, and consider the probability density functions $f_{\nu_1}$ and $f_{\nu_2}$, as defined by \eqref{eq_t_pdf}.  To verify that each $f_{\nu_j} \in L^2(\R)$, we simply note that 
$$
[f_{\nu_j}(x)]^2 = \dfrac{(1+x^2)^{-(2\nu_j+1)}}{[B(\nu_j,\tfrac12)]^2} = \dfrac{B(2\nu_j+\tfrac12,\tfrac12)}{[B(\nu_j,\tfrac12)]^2} f_{2\nu_j+\tfrac12}(x);
$$
this shows that $f_{\nu_j}^2 \in L^1(\R)$, equivalently, $f_{\nu_j} \in L^2(\R)$.  

Consequently $\widehat{f}_{\nu_j} \in L^2(\R)$ so, by Parseval's identity, 
\begin{align*}
\frac{1}{2\pi} \int_{-\infty}^\infty \dfrac{|t|^{\nu_1} K_{\nu_1}(|t|)}{2^{\nu_1-1} \Gamma(\nu_1)}
\dfrac{|t|^{\nu_2} K_{\nu_2}(|t|)}{2^{\nu_2-1} \Gamma(\nu_2)} \dd t 
&= \frac{1}{2\pi} \langle \widehat{f}_{\nu_1},\widehat{f}_{\nu_2} \rangle = \langle f_{\nu_1},f_{\nu_2} \rangle \\
&= \dfrac{1}{B(\nu_1,\tfrac12) B(\nu_2,\tfrac12)} \int_{-\infty}^\infty (1+x^2)^{-(\nu_1+\nu_2+1)} \dd x \\
&= \dfrac{B(\nu_1+\nu_2+\tfrac12,\tfrac12)}{B(\nu_1,\tfrac12) B(\nu_2,\tfrac12)}.
\end{align*}
In summary, we obtain the integral, 
\begin{align*}
\int_0^\infty t^{\nu_1+\nu_2} K_{\nu_1}(t) K_{\nu_2}(t) \dd t &= \pi \cdot \frac{1}{2\pi} \int_{-\infty}^\infty |t|^{\nu_1+\nu_2} K_{\nu_1}(|t|) K_{\nu_2}(|t|) \dd t \\
&= 2^{\nu_1+\nu_2-2} \pi \Gamma(\nu_1) \Gamma(\nu_2)\dfrac{B(\nu_1+\nu_2+\tfrac12,\tfrac12)}{B(\nu_1,\tfrac12) B(\nu_2,\tfrac12)} \\
&= 2^{\nu_1+\nu_2-2} \pi^{1/2} \frac{\Gamma(\nu_1+\nu_2+\tfrac12)}{\Gamma(\nu_1+\nu_2+1)} \Gamma(\nu_1+\tfrac12) \Gamma(\nu_2+\tfrac12),
\end{align*}
and, by analytic continuation, this result remains valid for all $\nu_1,\nu_2 \in \C$ such that $\Re(\nu_1) > -\tfrac12$, $\Re(\nu_2) > -\tfrac12$, and $\Re(\nu_1+\nu_2) > -1$.  We also note that the final result agrees with Gradshteyn and Ryzhik \cite[p.~660, 6.511(13); p.~684, 6.576(4)]{Gradshteyn}.
}\end{example}

The methods used above can be extended to calculate more general integrals of the form $\int_0^\infty t^\lambda K_{\nu_1}(at) K_{\nu_2}(bt) \dd t$.  This leads to final expressions involving the Gaussian hypergeometric function.

In closing this section, we consider probability density functions that are proportional to $(1-x^2)^{\nu_j-1}$, $-1 < x < 1$, $\nu_j > 0$, $j=1,2$.  Then we recover integrals of the type, $\int_0^\infty t^a J_{\nu_1}(t) J_{\nu_2}(t) \dd t$, the well-known integrals of Weber-Schafheitlin type \cite[p.~398]{Watson}.  

\medskip

\begin{example}
\label{ex_sphere}
{\rm 
For $\nu > \tfrac12$, let 
$$
f_\nu(x) = (1-x^2)^{\nu-\tfrac12},
$$
$-1 < x < 1$.  Then 
$$
\fhat_\nu(t) = \pi^{1/2} \Gamma(\nu+\tfrac12) (|t|/2)^{-\nu} J_\nu(|t|),
$$
$t \neq 0$.  
For $\nu_1, \nu_2 > \tfrac12$ we apply Parseval's identity to $f_{\nu_1}$ and $f_{\nu_2}$, obtaining 
\begin{align*}
\frac{1}{2\pi} \int_{-\infty}^\infty \pi^{1/2} & \Gamma(\nu_1+\tfrac12) (|t|/2)^{-\nu_1} J_{\nu_1}(|t|) \cdot \pi^{1/2} \Gamma(\nu_2+\tfrac12) (|t|/2)^{-\nu_2} J_{\nu_2}(|t|) \dd t \\
&= \frac{1}{2\pi} \langle \widehat{f}_{\nu_1},\widehat{f}_{\nu_2} \rangle = \langle f_{\nu_1},f_{\nu_2} \rangle \\
&= \int_{-1}^1 (1-x^2)^{\nu_1+\nu_2-1} \dd x 
= 2^{2(\nu_1+\nu_2)-1} B(\nu_1+\nu_2,\nu_1+\nu_2).
\end{align*}
Equivalently, 
\begin{multline*}
\Gamma(\nu_1+\tfrac12) \Gamma(\nu_2+\tfrac12) \frac{1}{2} \int_{-\infty}^\infty (|t|/2)^{-\nu_1} J_{\nu_1}(|t|) \cdot  (|t|/2)^{-\nu_2} J_{\nu_2}(|t|) \dd t \\
= 2^{2(\nu_1+\nu_2)-1} B(\nu_1+\nu_2,\nu_1+\nu_2) = 2^{2(\nu_1+\nu_2)-1} \frac{\Gamma(\nu_1+\nu_2)\Gamma(\nu_1+\nu_2)}{\Gamma(2(\nu_1+\nu_2))},
\end{multline*}
so 
\begin{align*}
\int_0^\infty t^{-\nu_1-\nu_2} J_{\nu_1}(t) J_{\nu_2}(t) \dd t 
&= 2^{\nu_1+\nu_2-1} \frac{\Gamma(\nu_1+\nu_2)\Gamma(\nu_1+\nu_2)}{\Gamma(2(\nu_1+\nu_2))\Gamma(\nu_1+\tfrac12) \Gamma(\nu_2+\tfrac12)} \\
&= \frac{\pi^{1/2} \Gamma(\nu_1+\nu_2)}{2^{\nu_1+\nu_2} \Gamma(\nu_1+\nu_2+\tfrac12) \Gamma(\nu_1+\tfrac12) \Gamma(\nu_2+\tfrac12)},
\end{align*}
where the latter equality follows by applying the duplication formula for the gamma function.  Further, by analytic continuation, the final result remains valid for  $\nu_1,\nu_2 \in \C$ such that $\Re(\nu_1+\nu_2) > 0$, and it agrees with Gradshteyn and Ryzhik \cite[p.~683, 6.575(2)]{Gradshteyn}.
}\end{example}

\section{Multidimensional Cauchy-beta integrals}
\label{sec_multidim_int}

Denote by $\calSd$ the vector space of $d \times d$ real symmetric matrices, and let $\calSdpd$ be the cone in $\calSd$ of positive definite matrices.  For $X \in \calSd$ denote by $\tr(X)$ and $\det X$ the trace and determinant, respectively, and also denote $\exp(\tr X)$ by $\etr(X)$.  For $\alpha \in \C$ such that $\Re(\alpha) > (d-1)/2$, define the gamma function for the cone $\calSdpd$
\begin{equation}
\label{eq_cone_gamma_funct}
\Gamma_d(\alpha) = \int_{\calSdpd} (\det X)^{\alpha} \etr(-X) \, (\det X)^{-(d+1)/2} \dd X,
\end{equation}
where $\dd X$ denotes Lebesgue measure on $\calSdpd$.  A well-known result, due to Wishart, and later to Ingham and Siegel, is that 
\begin{equation}
\label{eq_cone_gamma_funct_2}
\Gamma_d(\alpha) = \pi^{d(d+1)/4} \prod_{j=1}^d \Gamma(\alpha-\tfrac12(j-1)).
\end{equation}
Wishart's integral is now prominent in aspects of multivariate statistical analysis (see Muirhead \cite{Muirhead}) and has been generalized to the setting of the symmetric cones (see Ding, \textit{et al.} \cite{DGR}).  

For $\Sigma \in \calSdpd$ and $\alpha > (d-1)/2$, the function 
\begin{equation}
\label{eq_Wishart_pdf}
f_{\alpha,\Sigma}(X) = (\det \Sigma)^{-\alpha} (\det X)^{\alpha-\tfrac12(d+1)} \etr(-\Sigma^{-1} X)/\Gamma_d(\alpha),
\end{equation}
$X \in \calSdpd$, is a probability density function, and the corresponding probability distribution is called the \textit{Wishart distribution on} $\calSdpd$.  To prove that the function $f_{\alpha,\Sigma}$ integrates to $1$, i.e., that 
\begin{equation}
\label{eq_Wishart_int}
\int_\calSdpd (\det X)^\alpha \etr(-\Sigma^{-1} X) \, (\det X)^{-(d+1)/2} \dd X = \Gamma_d(\alpha) (\det \Sigma)^\alpha,
\end{equation}
one makes the change-of-variables $X \mapsto \Sigma^{1/2}X\Sigma^{1/2}$ in the integral, where $\Sigma^{1/2}$ denotes the unique positive definite square-root of $\Sigma$, apply the invariance of the measure $(\det X)^{-(d+1)/2} \dd X$ under such a substitution, and use the definition \eqref{eq_cone_gamma_funct} of $\Gamma_d(\alpha)$.  

Given a function $f \in L^1(\calSdpd)$, define the Fourier transform of $f$ to be 
$$
\fhat(T) = \int_\calSdpd \etr(-iTX) f(X) \dd X,
$$
$T \in \calSd$.  Applying classical Fourier transform theory, it follows that $\fhat \in L^2(\calSd)$ if $f \in L^2(\calSdpd)$; moreover, for $f_1, f_2 \in L^1(\calSdpd) \cap L^2(\calSdpd)$, the Parseval identity implies that 
\begin{multline}
\label{eq_matrix_Parseval}
\frac{1}{(2\pi)^{d(d+1)/2}} \int_\calSd \fonehat(T) \overline{\ftwohat(T)} \dd T \\
= \frac{1}{(2\pi)^{d(d+1)/2}} \langle \fonehat,\ftwohat \rangle 
= \langle f_1,f_2 \rangle = \int_\calSdpd f_1(X) \overline{f_2(X)} \dd X.
\end{multline}

Our first example of a Cauchy beta integral on the cone $\calSdpd$ extends \eqref{eq_Cauchy_beta_integral}.  In deriving this result, we denote by $\calS_d+i\calS_d = \{T_1 + iT_2: T_1,T_2 \in \calSd\}$ the vector space of complex symmetric $d \times d$ matrices.  

\medskip

\begin{example}
\label{example_matrix_Cauchy_beta}
\textit{Suppose that $\alpha_1,\alpha_2 \in \C$ such that $\Re(\alpha_1+\alpha_2) > d$, and let $\Sigma_1, \Sigma_2 \in \calS_d+i\calS_d$ such that $\Re(\Sigma_1), \Re(\Sigma_2) \in \calSdpd$.  Then 
\begin{multline}
\label{eq_matrix_Cauchy_beta}
\int_\calSd \det(I_d + i\Sigma_1 T)^{-\alpha_1} \det(I_d - i\Sigma_2 T)^{-\alpha_2} \dd T \\
= (2\pi)^{d(d+1)/2} \frac{\Gamma_d(\alpha_1+\alpha_2-\tfrac12(d+1))}{\Gamma_d(\alpha_1) \Gamma_d(\alpha_2)} \frac{(\det \Sigma_1)^{\alpha_2 - \tfrac12(d+1)} (\det \Sigma_2)^{\alpha_1 - \tfrac12(d+1)}}{(\det(\Sigma_1+\Sigma_2))^{\alpha_1+\alpha_2-\tfrac12(d+1)}}
\end{multline}
and 
\begin{equation}
\label{eq_matrix_Cauchy_beta2}
\int_\calSd \det(I_d + i\Sigma_1 T)^{-\alpha_1} \det(I_d + i\Sigma_2 T)^{-\alpha_2} \dd T = 0.
\end{equation}
}\end{example}

\Proof
For $\Sigma \in \calS_d+i\calS_d$, the right-hand side of \eqref{eq_Wishart_int} is analytic in $\Sigma$ for $\Re(\Sigma) \in \calSdpd$.  Therefore \eqref{eq_Wishart_int} remains valid for the larger region $\alpha > (d-1)/2$ and $\Re(\Sigma) \in \calSdpd$.  Consequently, it is a straightforward consequence of \eqref{eq_Wishart_int} that the Fourier transform of $f_{\alpha,\Sigma}$ is 
$$
\fhat_{\alpha,\Sigma}(T) = \det(I_d + i\Sigma T)^{-\alpha},
$$
$T \in \calSd$.  Further, by the same argument that proves $f_{\alpha,\Sigma} \in L^1(\calSdpd)$, we find that if $\alpha > d/2$ and $\Re(\Sigma) \in \calSdpd$ then $f_{\alpha,\Sigma} \in L^2(\calSdpd)$, in which case we also obtain $\fhat_{\alpha,\Sigma} \in L^2(\calSd)$.

Now suppose that $\alpha_1,\alpha_2 \in \R$ with $\alpha_1, \alpha_2 > d/2$.  For $\Sigma_1, \Sigma_2 \in \calSdpd$, we apply the Parseval identity \eqref{eq_matrix_Parseval} to obtain 
\begin{align*}
\frac{1}{(2\pi)^{d(d+1)/2}}& \int_\calSd \det(I_d + i\Sigma_1 T)^{-\alpha_1} \det(I_d - i\Sigma_2 T)^{-\alpha_2} \dd T \\
&= \frac{1}{(2\pi)^{d(d+1)/2}} \langle \widehat{f}_{\alpha_1},\widehat{f}_{\alpha_2} \rangle = \langle f_{\alpha_1},f_{\alpha_2} \rangle \\
&= \frac{(\det \Sigma_1)^{-\alpha_1} (\det \Sigma_2)^{-\alpha_2}}{\Gamma_d(\alpha_1) \Gamma_d(\alpha_2)} \int_\calSdpd (\det X)^{\alpha+\alpha_2-d-1} \etr(-(\Sigma_1^{-1}+\Sigma_2^{-1}) X) \dd X.
\end{align*}
On applying \eqref{eq_Wishart_int} to evaluate the latter integral, and then simplifying the resulting expression, we obtain \eqref{eq_matrix_Cauchy_beta}.  

Next, we apply the analyticity of the right-hand side of \eqref{eq_matrix_Cauchy_beta} as a function of $\Sigma_1, \Sigma_2$ to deduce that the result remains valid for $\Sigma_1, \Sigma_2 \in \calS_d+i\calS_d$ such that $\Re(\Sigma_1), \Re(\Sigma_2) \in \calSdpd$.  

For complex $\alpha_1, \alpha_2$, it follows from \eqref{eq_cone_gamma_funct_2} that the right-hand side of \eqref{eq_Cauchy_beta_integral} is analytic in $\alpha_1, \alpha_2$ for $\Re(\alpha_1+\alpha_2) > d$.  Therefore \eqref{eq_Cauchy_beta_integral} remains valid on the stated region, i.e., for $\alpha_1,\alpha_2 \in \C$ such that $\Re(\alpha_1+\alpha_2) > d$, and $\Sigma_1, \Sigma_2 \in \calS_d+i\calS_d$ such that $\Re(\Sigma_1), \Re(\Sigma_2) \in \calSdpd$.  

As for \eqref{eq_matrix_Cauchy_beta2}, that result follows by an argument entirely analogous to the derivation of \eqref{eq_Cauchy_beta_integral2}.  
$\qed$

\medskip

The next result, a beta integral of Cauchy--Selberg type, is a special case of an integral due to Selberg.  See Morris \cite[p.~92]{Morris} for more general instances of these integrals, an account of some of their history, and applications to constant term identities.  
 
\medskip

\begin{example} (Selberg.)
\textit{Suppose that $\alpha_1,\alpha_2,\sigma_1,\sigma_2 \in \C$ such that $\Re(\alpha_1+\alpha_2) > d$ and $\Re(\sigma_1), \Re(\sigma_2) > 0$.  Then 
\begin{multline}
\label{eq_eigens_Cauchy_beta}
\int_{\R^d} \prod_{1 \le j < k \le d} |t_j - t_k| \cdot \prod_{j=1}^d (1+i\sigma_1 t_j)^{-\alpha_1} (1 - i\sigma_2 t_j)^{-\alpha_2} \dd t_j \\
= d! \, 2^{d(d+1)/2} \pi^{d/2} \frac{\Gamma_d(\tfrac12d) \Gamma_d(\alpha_1+\alpha_2-\tfrac12(d+1))}{\Gamma_d(\alpha_1) \Gamma_d(\alpha_2)} \frac{\sigma_1^{\alpha_2 d - \tfrac12d(d+1)} \sigma_2^{\alpha_1 d - \tfrac12d(d+1)}}{(\sigma_1+\sigma_2)^{(\alpha_1+\alpha_2)d - \tfrac12d(d+1)}}
\end{multline}
and 
\begin{equation}
\label{eq_eigens_Cauchy_beta1}
\int_{\R^d} \prod_{1 \le j < k \le d} |t_j - t_k| \cdot \prod_{j=1}^d (1+i\sigma_1 t_j)^{-\alpha_1} (1 + i\sigma_2 t_j)^{-\alpha_2} \dd t_j = 0.
\end{equation}
}\end{example}

\Proof
Suppose that $\alpha_1,\alpha_2 \in \R$, where $\alpha_1, \alpha_2 > (d-1)/2$.  Also suppose that $\sigma_1, \sigma_2 \in \C$ such that $\Re(\sigma_1), \Re(\sigma_2) > 0$.  Letting $I_d$ denote the $d \times d$ identity matrix, we return to \eqref{eq_matrix_Cauchy_beta} and substitute $\Sigma_1 = \sigma_1 I_d$ and $\Sigma_2 = \sigma_2 I_d$.  Then we obtain 
\begin{multline}
\label{eq_eigens_Cauchy_beta2}
(2\pi)^{d(d+1)/2} \frac{\Gamma_d(\alpha_1+\alpha_2-\tfrac12(d+1))}{\Gamma_d(\alpha_1) \Gamma_d(\alpha_2)} \frac{\sigma_1^{\alpha_2 d - \tfrac12d(d+1)} \sigma_2^{\alpha_1 d - \tfrac12d(d+1)}}{(\sigma_1+\sigma_2)^{(\alpha_1+\alpha_2)d - \tfrac12d(d+1)}} \\
= \int_\calSd \det(I_d + i\sigma_1 T)^{-\alpha_1} \det(I_d - i\sigma_2 T)^{-\alpha_2} \dd T.
\end{multline}
Let $O(d)$ denote the group of $d \times d$ orthogonal matrices.  For generic $H \in O(d)$, denote by $\dd H$ the corresponding Haar measure, normalized to have total volume $1$.  Each $T \in \calSd$ can be expressed in the form 
\begin{equation}
\label{eq_eigens_Cauchy_beta3}
T = H \, \diag(t_1,\ldots,t_d) \, H'
\end{equation}
where $\diag(t_1,\ldots,t_d)$ is a diagonal matrix whose diagonal entries are the eigenvalues of $T$.  By Muirhead \cite[p.~103]{Muirhead}, the Lebesgue measure on $\calSd$ decomposes as 
\begin{equation}
\label{eq_eigens_Cauchy_beta4}
\dd T = \frac{\pi^{d^2/2}}{d! \, \Gamma_d(\tfrac12d)} \, \dd H \prod_{1 \le j < k \le d} |t_j - t_k| \, \dd t_1 \cdots \dd t_d.
\end{equation}
Applying \eqref{eq_eigens_Cauchy_beta3} and \eqref{eq_eigens_Cauchy_beta4} to the integrand in \eqref{eq_eigens_Cauchy_beta2}, and integrating over $H$ using the fact that the Haar measure has volume $1$, we find that the integral in \eqref{eq_eigens_Cauchy_beta2} equals  
$$ 
\frac{\pi^{d^2/2}}{d! \, \Gamma_d(\tfrac12d)} \int_{\R^d} \prod_{1 \le j < k \le d} |t_j - t_k| \cdot \prod_{j=1}^d (1+i\sigma_1 t_j)^{-\alpha_1} (1 - i\sigma_2 t_j)^{-\alpha_2} \dd t_j.
$$
On equating the latter expression with the left-hand side of \eqref{eq_eigens_Cauchy_beta2} and simplifying the normalizing constant, we obtain \eqref{eq_eigens_Cauchy_beta}.  

To complete the proof, we again observe from \eqref{eq_cone_gamma_funct_2} that the function $\Gamma_d(\alpha)$, $\alpha \in \C$, is analytic on the region $\Re(\alpha) > (d-1)/2$; then it follows that the right-hand side of \eqref{eq_eigens_Cauchy_beta} is analytic in $\alpha_1, \alpha_2$ for $\Re(\alpha_1+\alpha_2) > d-1$.  

Finally, \eqref{eq_eigens_Cauchy_beta1} is proved by starting with \eqref{eq_matrix_Cauchy_beta2} and proceeding analogously to the proof of \eqref{eq_eigens_Cauchy_beta}.
$\qed$

\section{Concluding remarks}
\label{sec_conclusions}

The results in Sections \ref{sec_Cauchy} and \ref{sec_multidim_int} can be extended in numerous directions, and at the heart of each such extensions lies an analog of the Parseval identity.  For instance, $q$-series variations on the results of Section \ref{sec_Cauchy} can be obtained using $q$-analogs of the Fourier transform and their corresponding Parseval identities; see, e.g., Rubin \cite{Rubin}.  

The analysis in Section \ref{sec_multidim_int} on the cone of positive definite matrices can be extended \textit{via} the results of Gindikin \cite{Gindikin} to some general homogeneous domains.  As a consequence, we would obtain generalizations of Example \ref{example_matrix_Cauchy_beta} in which the integrands extending \eqref{eq_matrix_Cauchy_beta} involve products of powers of principal minors of the matrices $I_d + i \Sigma$ and $I_d - i\Sigma_2 T$.  On specializing such extensions to the case in which the domains of integration are the symmetric cones \cite{DGR}, we will obtain generalizations of Example \ref{ex_Cauchy_beta} in which the integrand contains zonal (spherical) polynomials.  Further, we will obtain generalizations of Examples \ref{example_1F1}--\ref{ex_sphere} involving the spherical functions, the Bessel functions, or the confluent hypergeometric functions defined on the symmetric cones.  

In other directions, analogs of the results given here can be obtained using other integral transforms, such as the Mellin and Hankel transforms, then applying Parseval-Plancherel identities for those transforms.  We have already noted that Askey and Roy \cite{AskeyRoy} and Marichev \cite{Marichev} derived several integral formulas using the Parseval identity for the classical Mellin transform, and the possibility exists that the Parseval identity for the classical Hankel transform can be applied similarly.  Moreover these results suggest that numerous multivariate integrals over matrix spaces can be evaluated using multidimensional Parseval identities for the Hankel transforms of functions of matrix argument and for the spherical Mellin transforms on the symmetric cones \cite{DGR}.  

In conclusion, we see these results as evidence of Dick Askey's foresight and vision; he saw in the seemingly simple Cauchy beta integral an enormous breadth and depth of mathematics.


\label{biblio}
\bibliographystyle{ims}

\end{document}